\documentstyle[11pt]{article}
\newtheorem{Lemma}{Lemma}
\newtheorem{Proposition}[Lemma]{Proposition}

\newtheorem{Theorem}[Lemma]{Theorem}
\newtheorem{Conjecture}[Lemma]{Conjecture}

\newtheorem{Corollary}[Lemma]{Corollary}

\newcommand{\EE}{\mbox{${\cal E}$}}

\newcommand{\bX}{{\bf X}}

\newcommand{\bx}{{\bf x}}

\newcommand{\sfrac}[2]{{\textstyle\frac{#1}{#2}}}
\newcommand{\qed}{\ \ \rule{1ex}{1ex}}

\newcommand{\JJ}{{\cal J}}
\newcommand{\XX}{{\cal X}}
\newcommand{\bW}{{\bf W}}
\newcommand{\bw}{{\bf w}}
\newfam\Bbbfam
\font\tenBbb=msbm10
\font\sevenBbb=msbm7
\font\fiveBbb=msbm5
\textfont\Bbbfam=\tenBbb
\scriptfont\Bbbfam=\sevenBbb
\scriptscriptfont\Bbbfam=\fiveBbb
\def\Bbb{\fam\Bbbfam \tenBbb}
\newcommand{\Zbold} {{ {\Bbb Z} }}

\begin{document}
\title{The Asymmetric One-Dimensional Constrained Ising Model:
Rigorous Results}

\author{David Aldous \thanks{Department of Statistics,
	367 Evans Hall,
	University of California,
	Berkeley CA 94720}
\and Persi Diaconis \thanks{Department of Statistics,
Sequoia Hall,
Stanford University CA 94305}
}

\maketitle

{\em Running head}: Constrained Ising Model.

\vspace{1.0in}

{\em Designated author}

David Aldous

Department of Statistics

367 Evans Hall \# 3860

U.C. Berkeley CA 94720-3860

phone: 510-642-2781

FAX: 510-642-7892

e-mail: aldous@stat.berkeley.edu

\newpage
\begin{abstract}
We study a one-dimensional spin (interacting particle) system, with
product Bernoulli($p$) stationary distribution, in which a site can
flip only when its left neighbor is in state $+1$.
Such models have been studied in physics as simple exemplars of systems
exhibiting slow relaxation.
In our ``East" model the natural conjecture is that the relaxation time
$\tau(p)$, that is 1/(spectral gap), satisfies
$\log \tau(p) \sim \frac{\log^2 (1/p)}{\log 2}
\mbox{ as } p \downarrow 0$.
We prove this up to a factor of $2$.
The upper bound uses the Poincar{\'e} comparison argument applied to a 
``wave" (long-range) comparison process, which we analyze by probabilistic
techniques.  
Such comparison arguments go back to Holley (1984, 1985).
The lower bound, which atypically is not easy, involves construction and analysis of 
a certain ``coalescing random jumps" process.
\end{abstract}

\vspace{0.3in}


\newpage
\section{Introduction}
\label{INT}

The asymmetric one-dimensional constrained Ising model,
or more briefly the {\em East process},
is an interacting particle system with sites $\Zbold^1$
and each site having two states $\{0,1\}
= \{ {\rm unoccupied, occupied}\}$.
Its essential qualitative feature is that
a site can change state only when the site to its left
is occupied.
The flip rates at each site $i$ are specified by:

\noindent
$\bullet$ if site $i-1$ is in state $0$ then state $i$ cannot change;

\noindent
$\bullet$ if site $i-1$ is in state $1$ then state $i$ flips
$\begin{array}{cc}
0 \to 1; &\mbox{  rate } p\\
1 \to 0; & \mbox{ rate }1-p
\end{array}$\\
where $0<p<1$ is a parameter.
Here is an equivalent description.
Each particle, at rate $1$
(that is, the times of a Poisson$(1)$ process)
sends a ``pulse" to the site to its East, and the state of that
site is reset via
\[ P({\rm occupied}) = p, \quad
P({\rm unoccupied}) = 1-p. \]
A careful construction of the process is outlined in
the Appendix.  
By routine arguments the i.i.d. Bernoulli$(p)$ measure is
the unique non-trivial stationary distribution and this {\em East process} is reversible.
Write $(\bX(t), 0 \leq t < \infty) = ((X_i(t), i \in \Zbold), 0 \leq t < \infty)$
for the stationary process.
We are interested in studying the {\em relaxation time}
(defined as 1/spectral gap),
say $\tau(p)$,
as $p \downarrow 0$.
This specific process, and questions concerning its relaxation time, were introduced by J\"{a}ckle and Eisinger
\cite{JE91}
and further studied in the physics literature in
\cite{EJ93,mauch,MGIP}.
Sollich and Evans \cite{sollich} argue non-rigorously that
\begin{equation}
\log \tau(p) \sim \frac{\log^2 (1/p)}{\log 2}
\mbox{ as } p \downarrow 0
\label{conj}
\end{equation}
({\em log} is natural logarithm),
and observe that (\ref{conj}) is consistent with Monte Carlo simulations.
Section \ref{sec-HA} outlines a simple heuristic argument.
The purpose of this paper is prove rigorous bounds  
which support this conjecture.
\begin{Theorem}
\label{T1}
(a) $\log \tau(p) \leq \left(\sfrac{1}{\log 2} + o(1)
\right) \log^2 (1/p)
\mbox{ as } p \downarrow 0 $.\\
(b) $\log \tau(p) \geq \left(\sfrac{1}{2 \log 2} - o(1)
\right) \log^2 (1/p)
\mbox{ as } p \downarrow 0 $.
\end{Theorem}
The proof of the lower bound (b), outlined in section 
\ref{sec-overview} with details in section \ref{sec-pfL},
involves the usual method of applying the variational
characterization (\ref{EC}) to a suitable test function $g$.
But unusually, finding a good $g$ is not intuitively simple,
and in fact we define $g$ only implicitly in terms of a certain
{\em coalescing random jumps} process which we need to invent.
Section \ref{sec-remarks} motivates this particular process.
A more inspired choice of $g$ might allow one
to remove the factor $\sfrac{1}{2}$ from the lower bound and thus
prove the conjecture (\ref{conj}).

The proof of the upper bound (a) uses the
Poincar{\'e} comparison method, 
used by Holley \cite{holley82,holley85} and developed by
Diaconis and Saloff-Coste \cite{DSC93a},
which bounds the relaxation time of one process
in terms of the relaxation time of another ``comparison" 
process with the same stationary distribution.
Though the simplest comparison process would 
be the process in which sites
flip independently, 
making the comparison with this process seems technically difficult.
Instead, following an idea used by Holley \cite{holley85}, we use for comparison a certain ``long range" process, which for us is (section \ref{sec-wave}) the following {\em wave} process,
For each particle (at site $i$, say) at rate $1$ a wave appears
which, in the $10/p$ sites to the right of $i$,
deletes existing particles and replaces them by an independent
Bernoulli($p$) family of particles.
Proposition \ref{Pwave} shows that for small $p$ the relaxation time
of the wave process is bounded by the constant $10/3$.
The proof (section \ref{sec-mg}) uses probabilistic methods:
coupling and supermartingale estimates for the position
of the rightmost particle in the 
wave process started from a single particle.
Section \ref{sec-comp}
gives the argument comparing the East process with the wave process.

\subsection{Remarks on related models and techniques}
{\bf 1.} 
The East model is a special case of more general {\em constrained}
or {\em facilitated} spin models on lattices, in which a site is
permitted to flip only when some prescribed number of neighbors are in
state $1$.  Such models go back to Frederickson and Anderson \cite{FA84},
and a recent review is in Pitts et al \cite{PYA00}.  The models are intended
to illustrate the slow relaxation behavior of liquids near
the glass transition.  Our methods are tied closely
to the specifics of the East model, but making rigorous the treatments
in \cite{PYA00} concerning other models would be an interesting challenge.
A specific next model one might study is to take as site space the
infinite rooted binary tree, and allow flips at a site only when
both children of that site are in state $1$.
In this model it is easy to show
\[ \tau(p) \geq c (p - \sfrac{1}{2})^{-2}, \ p > 1/2 \]
and one would like to prove the corresponding upper bound.
Apparently much harder is the {\em North-East model}, 
where site space is $\Zbold^2$ and flips are allowed at $(x,y)$
only if sites $(x-1,y)$ and $(x,y-1)$ are both in state $1$.
Here one would like to 
find a critical exponent $\alpha > 0$ such that
\[ \tau(p) \sim c (p - p_{{\rm crit}})^{-\alpha} \mbox{ as }
p \downarrow p_{{\rm crit}} \]
where
$p_{{\rm crit}}$ is the critical value for oriented percolation
of $0$-sites.

In the terminology of \cite{PYA00} our model should be called the {\em West} model; 
we reversed direction
to make site labeling more natural.

{\bf 2.} 
Since Diaconis and Stroock \cite{DS90} there have been many papers
studying relaxation times for Markov chains on various
``combinatorial" state spaces via the basic comparison method
(implicitly comparing the given chain with the i.i.d. chain).
Comparisons with some judiciously-chosen alternate chain were
used in \cite{DSC93a} to study the exclusion process,
and in \cite{RT-glauber} for other statistical physics models.
Martinelli \cite{martinelli94} uses such arguments in studying
the subcritical two-dimensional Ising process.

{\bf 3.}
One might suppose that general results for one-dimensional Ising-type
models implied \begin{equation}\tau (p) < \infty \mbox{ for each } p > 0
\label{taufinite}
\end{equation}
But the usual
general results such as Holley \cite{holley85} 
implicitly assume flip rates are strictly positive and so do not
apply to the East process.
So even (\ref{taufinite}) may be a new result.
Note also that the East process is not {\em monotone} (or
{\em attractive}) in the sense of interacting particle systems
\cite{lig85}.  Our {\em extension coupling} in section
\ref{sec-CSP} uses a weaker structural property.

{\bf 4.}
For any stationary reversible Markov process
${\bX}(t)$, there are two useful ways of viewing the relaxation time $\tau$.
One is via the ``infinitesimal time" 
Rayleigh-Ritz variational characterization
(e.g. \cite{bremaud99} eq. 6.2.10):
\begin{equation}
\tau = \sup_g \frac{{\rm var\ } g}{{\cal E}(g,g)}
\label{EC}
\end{equation}
where ${\rm var\ } g $ is the variance of $g(\bX(0))$; and
\[ \EE(g,g) = \sfrac{1}{2} \lim_{t \downarrow 0}
t^{-1} E(g(\bX(t))-g(\bX(0)))^2 . \]
The other is via the finite-time maximal correlation property:
\begin{equation}
\sup_{f,g} {\rm cor}(f(\bX(0)),g(\bX(t)) = \exp(-t/\tau), \quad 
0 \leq t < \infty . \label{max-cor}
\end{equation}
In general the {\em sup}s are over all square-integrable functions;
in our setting of interacting particle systems, in (\ref{max-cor}) we may restrict to
bounded functions depending on only finitely many sites.
See the Appendix.
Our arguments are purely finitistic -- we derive estimates of the
spectral gap for the East process on states 
$\{0,1,\ldots,n\}$ which are uniform in large $n$ -- and we downplay
the routine abstract arguments needed to pass to the limit infinite-site
models.

{\bf 5.}
As well as the (mathematically natural) relaxation time $\tau(p)$
of the whole East process,
it is physically natural to consider the single-spin asymptotic
relaxation time
$\tau^0(p)$ defined by
\[ P(X_0(t) = 1|X_0(0) = 1) - p  
= \exp(-\sfrac{t}{\tau^0(p) \pm o(1)})
\mbox{ as } t \to \infty . \]
One immediately has
$\tau^0(p) \leq \tau(p)$
and one expects equality, but there is no simple general proof
of equality.

{\bf 6.}
The procedure for deriving the East process from Glauber dynamics
for product measure has a natural abstraction to arbitrary
reversible Markov processes.
Taking for simplicity finite state space ${\cal X}$,
let $Q(x,y)$ be the transition rate matrix and $\pi$ the stationary
distribution of a reversible chain, and let $S$
be a symmetric subset of ${\cal X} \times {\cal X}$.
Then the derived chain
\[
\hat{Q}(x,y) = Q(x,y)  
1_{(y \neq x, \ (x,y) \in S)} \]
is reversible with the same stationary distribution $\pi$.
The variational principle immediately implies that the relaxation
time for the derived chain cannot be smaller than for the original chain.
Upper bounding such relaxation times, when the original chain
is e.g. the usual Glauber dynamics for the Ising model,
is a natural open problem.

\section{Heuristics}
\label{sec-heur}
Our proofs are indirect, because we are unable to do sharp calculations
directly with the East process.
Underlying our proofs is an intuitive visualization of how the
process evolves, which we shall outline in this section.

\subsection{Visualizing a realization}
Fix small $p$.
A typical particle at a typical time is isolated from other
particles by a distance of order $1/p$.
Figure 1 illustrates schematically a realization of
the East process
from an isolated particle at site $0$.

\setlength{\unitlength}{0.16in}
\begin{picture}(15,25)(-3,0)
\multiput(0,22)(0,-1){23}{\circle*{0.4}}
\multiput(1,21)(0,-2){11}{\circle*{0.4}}
\put(2,0){\circle*{0.4}}
\put(1,0){\circle*{0.4}}
\multiput(5,22)(0,-1){23}{\circle*{0.4}}
\put(6,22){\circle*{0.4}}
\multiput(7,22)(0,-1){5}{\circle*{0.4}}
\put(8,20){\circle*{0.4}}
\put(6,18){\circle*{0.4}}
\multiput(6,17)(0,-2){9}{\circle*{0.4}}
\put(7,0){\circle*{0.4}}
\put(6,0){\circle*{0.4}}
\put(9,12){more}
\put(8.9,11){cycles}
\put(9.5,20){\vector(0,-1){18}}
\multiput(12,22)(0,-1){23}{\circle*{0.4}}
\multiput(13,22)(0,-1){2}{\circle*{0.4}}
\multiput(13,20)(0,-4){6}{\circle*{0.4}}
\multiput(14,22)(0,-1){2}{\circle*{0.4}}
\multiput(15,21)(0,-1){22}{\circle*{0.4}}
\multiput(16,18)(0,-4){5}{\circle*{0.4}}
\multiput(20,22)(0,-1){23}{\circle*{0.4}}
\multiput(21,22)(0,-1){2}{\circle*{0.4}}
\multiput(22,21)(0,-1){6}{\circle*{0.4}}
\multiput(23,22)(0,-1){3}{\circle*{0.4}}
\put(23,18){\circle*{0.4}}
\multiput(21,15)(0,-2){8}{\circle*{0.4}}
\put(21,16){\circle*{0.4}}
\put(-0.23,23){0}
\put(0.77,23){1}
\put(1.77,23){2}
\put(4.77,23){0}
\put(5.77,23){1}
\put(6.77,23){2}
\put(7.77,23){3}
\put(11.77,23){0}
\put(12.77,23){1}
\put(13.77,23){2}
\put(14.77,23){3}
\put(15.77,23){4}
\put(19.77,23){0}
\put(20.77,23){1}
\put(21.77,23){2}
\put(22.77,23){3}
\put(-1,22){\vector(0,-1){6}}
\put(-4,20){time}
\end{picture}

\vspace{0.12in}

Fig. 1. 
{\small 
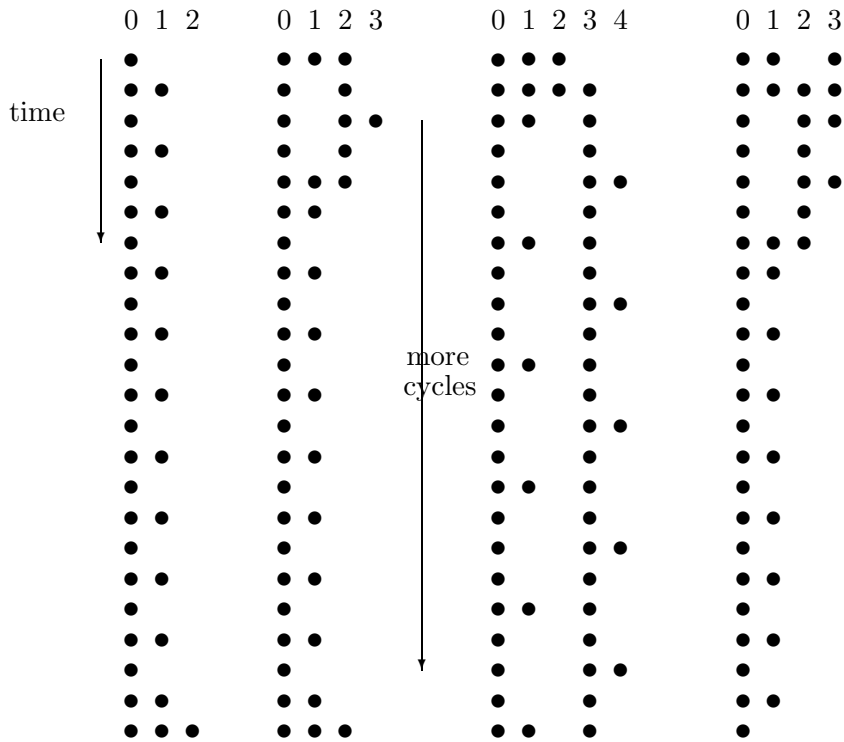
Schematic representation of successive configurations in the
East process started from a particle at site $0$.
Columns represent consecutive time intervals.
 From configuration $\bullet \bullet$ the process is much more likely
to transition to $\bullet$ than to $\bullet \bullet \bullet$;
column $1$ represents an order $1/p$ number of cycles
$\bullet \to \bullet \bullet \to \bullet$
before first attaining $\bullet \bullet \bullet$.
 From there, it is unlikely to transition to 
$\bullet \bullet \bullet \bullet$
but instead is equally likely to transition to
$\bullet \bullet$ or to $\bullet \circ \bullet$.
Column 2 indicates that either of the latter options will likely
cause the process to return to $\bullet$.
It therefore requires order $1/p$ more cycles
$\bullet \to \bullet  \bullet \bullet \to \bullet$
before first attaining $ \bullet \bullet \bullet \ \bullet$.
 From there it is possible to reach a locally recurrent 
configuration $\bullet \circ \circ \ \bullet$
(column 3) but the process will subsequently return to
$\bullet \bullet \bullet \ \bullet$
and the other likely possibilities lead (column 4) back to
$\bullet$.}

\vspace{0.09in}
\noindent
The moral of figure 1 is that, starting from an isolated particle, it takes a long time to reach a configuration
with another particle at some distance $1 \ll d \ll 1/p$
from the original particle;
and such configurations tend to relapse back to the original
configuration.  

\subsection{Minimum-energy paths}
\label{sec-ME}
For this section, take the site-space to be $\Zbold^+$ 
but consider only finite configurations.  Write
$|\bx| := \sum_i x_i$ for the
number of particles in configuration $\bx$.
A natural metric (in fact, an ultrametric) on configuration space 
is
\begin{equation}
h(\bx,\bx^\prime) := 
\min_{ \mbox{{\scriptsize paths $(\bx^j)$ from $\bx$ to $\bx^\prime$}}} \quad
 \max_j |\bx^j|
\label{hdef}
\end{equation}
where the minimum is over paths
$\bx = \bx^0,\bx^1,\bx^2,\ldots,\bx^\prime$
such that each successive pair
$(\bx^j,\bx^{j+1})$ is a possible transition of the East process.
This is natural because, in the usual statistical physics
picture, the relative probabilities of configurations $\bx$ 
depend only on $|\bx|$, which we can interpret as determining the
{\em energy} of $\bx$, so that $h(\bx,\bx^\prime)$ is the
(absolute) height of the ``energy barrier" separating the configurations.
Note that, in contrast to most interacting particle models,
if $\bx^\prime$ differs from $\bx$ only by virtue of moving an isolated
particle to an adjacent site, then $\bx$ and $\bx^\prime$
are {\em not} close in the $h$ metric; there is no short or low-energy path
between them.

Write $\delta_0 + \delta_i + \ldots$ for the configuration
consisting of particles at sites $0$ and $i$ and \ldots \ .
Figure 2 illustrates a recursive construction, which easily
proves
the following
(noted in the physics literature, e.g 
\cite{mauch} section 4).
Chung et al \cite{CDG01} give further analysis of such combinatorics.

\setlength{\unitlength}{0.12in}
\begin{picture}(-33,13)(-5,0)
\multiput(-5,9)(0,-1){10}{\circle*{0.4}}
\multiput(-4,8)(0,-1){2}{\circle*{0.4}}
\multiput(-4,2)(0,-1){2}{\circle*{0.4}}
\multiput(-3,7)(0,-1){6}{\circle*{0.4}}
\multiput(-2,5)(0,-1){2}{\circle*{0.4}}
\multiput(-1,4)(0,-1){5}{\circle*{0.4}}
\multiput(10,9)(0,-1){10}{\circle*{0.4}}
\multiput(14,9)(0,-1){10}{\circle*{0.4}}
\multiput(15,8)(0,-1){2}{\circle*{0.4}}
\multiput(15,2)(0,-1){2}{\circle*{0.4}}
\multiput(16,7)(0,-1){6}{\circle*{0.4}}
\multiput(17,5)(0,-1){2}{\circle*{0.4}}
\multiput(18,4)(0,-1){5}{\circle*{0.4}}
\multiput(25,9)(0,-1){10}{\circle*{0.4}}
\multiput(33,9)(0,-1){10}{\circle*{0.4}}
\multiput(26,1)(0,1){2}{\circle*{0.4}}
\multiput(26,7)(0,1){2}{\circle*{0.4}}
\multiput(27,2)(0,1){6}{\circle*{0.4}}
\multiput(28,4)(0,1){2}{\circle*{0.4}}
\multiput(29,5)(0,1){5}{\circle*{0.4}}
\put(-5.23,10){0}
\put(-4.23,10){1}
\put(-3.23,10){2}
\put(-2.23,10){3}
\put(-1.23,10){4}
\put(-0.23,10){5}
\put(0.77,10){6}
\put(1.77,10){7}
\put(2.77,10){8}
\put(9.77,10){0}
\put(10.77,10){1}
\put(11.77,10){2}
\put(12.77,10){3}
\put(13.77,10){4}
\put(14.77,10){5}
\put(15.77,10){6}
\put(16.77,10){7}
\put(17.77,10){8}
\put(24.77,10){0}
\put(25.77,10){1}
\put(26.77,10){2}
\put(27.77,10){3}
\put(28.77,10){4}
\put(29.77,10){5}
\put(30.77,10){6}
\put(31.77,10){7}
\put(32.77,10){8}
\end{picture}

\vspace{0.12in}

Fig. 2. 
{\small 
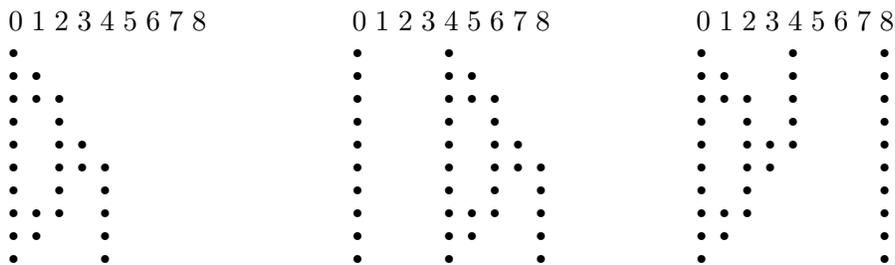
A path from $\delta_0$ to $\delta_0 + \delta_4$
(left column) is used to define 
a path from $\delta_0$ to $\delta_0 + \delta_8$.
The latter path has length $3$ times the former length,
and uses one extra particle.
}

\begin{Lemma}
\label{Lm}
Let $m \geq 0$.
There exists a path from
configuration $\delta_0$
to configuration $\delta_0 + \delta_{2^m}$
of length $3^m$
such that 
the maximum of $|\bx|$ over configurations $\bx$ on the path
equals $m+2$.
In particular
\[ h(\delta_0,\delta_0 + \delta_{2^m}) \leq m+2 . \]
\end{Lemma}

\subsection{Heuristic analysis of $\tau(p)$}
\label{sec-HA}
A very crude heuristic analysis (in part related to discussion
in \cite{sollich}) involves the following two ideas.
\\
(i) Write $s(p,j)$ for the mean time
for the East process started at $\delta_0$ to reach some
configuration with site $j$ occupied.
Since particles are typically separated by distance $1/p$,
the relaxation time should be roughly $s(p,1/p)$, since this
is the time until a particle has some influence on its initial
nearest neighbor particle.
\\
(ii) Lemma \ref{Lm} suggests that getting from configuration 
$\delta_0$ to a configuration with a particle at site $j$
involves an ``excursion" going over an energy barrier of height roughly
$h = \log_2 j$, in other words through configurations of
probability roughly $p^h$.
An excursion required to pass through states of some exponentially
small probability $q$ should take time roughly $1/q$.
So $s(p,j)$ should be
roughly $(1/p)^h$.

Combining (i) and (ii) suggests $\tau(p) \approx s(p,1/p) \approx 
(1/p)^{\log_2 1/p}$, which is (\ref{conj}).

\subsection{A continuum limit process?}
We would like to prove the following conjecture
(in which the $t(p)$ would be, up to constant factors, the
relaxation times $\tau(p)$)
concerning the extent of the ``excursions" above.
\begin{Conjecture}
Let $\bX^0(t)$ be the 
East process on sites $\Zbold^+$, started from a single particle at site $0$.
Let $r(\bx)$ be the rightmost occupied site in configuration $\bx$.
Then there exist constants $t(p)$ and a function $G(x), 0<x<\infty$ such that,
for all $\omega(p) \downarrow 0$ sufficiently slowly,
\begin{equation}
P \left( \sup_{0 \leq s \leq \omega(p)t(p)} r(\bX^0(s)) > \sfrac{x}{p} \right)
\sim G(x)\omega(p) 
\mbox{ as } p \downarrow 0 \mbox{ with $x$ fixed.}
\end{equation}
\end{Conjecture}
If this were true with suitable tail conditions on $G(\cdot)$, 
it would follow quite easily
(using the qualitative properties of section \ref{sec-CSP})
that as $p \downarrow 0$ the stationary East process on $\Zbold$, rescaled
by taking sites to be distance $p$ apart and speeding up time
by $t(p)$, 
converges to
a limit process
$(\Xi(t), 0 \leq t < \infty)$
described as follows.\\
(a) At fixed time $t$, $\Xi(t)$ is a Poisson (rate $1$) process of
particles on $R$.\\
(b) Each particle creates ``waves" of random lengths, the rate of
creation of waves of length $>l$ being $G(l)$.\\
(c) A wave $(x,x+l]$ instantaneously deletes all particles in $(x,x+l]$ and
replaces them by a Poisson (rate $1$) process of particles on 
$(x,x+l]$.

But proving the conjecture seems difficult; we do not have
even a heuristic derivation of a formula for $G(x)$.

\section{The wave process}
\label{sec-wave}
In this section we take the site space to be
$\Zbold^+ := \{0,1,2,\ldots\}$.
We generalize the East process by introducing, in addition to
$0<p<1$, a second parameter $v \in \{1,2,3,\ldots\}$,
and prescribing transition rates as follows.
\begin{quote}
Each particle (at site $i$, say)
at rate $1$ creates a ``wave" of length $v$, which instantaneously
deletes any particles in sites 
$\{i+1,i+2,\ldots,i+v\}$
and replaces them by an i.i.d. Bernoulli($p$) 
process of particles on sites
$\{i+1,i+2,\ldots,i+v\}$.
\end{quote}
Also specify that site $0$ is always occupied.
We call this the {\em wave process},
and write it as $\bW(t) = (W_i(t), \ i \in \Zbold^+)$.
The particular case $v=1$ is the East process.
It is easy to check that the general wave process is reversible
and its stationary distribution is
i.i.d. Bernoulli($p$) 
on sites $\{1,2,\ldots\}$ with the fixed particle at site $0$.

In section \ref{sec-CSP} we give some qualitative properties which do not
depend on $v$; then in section \ref{sec-mg} we give a quantitative bound
on the relaxation time for large $v$.

If $I$ is an interval of sites then we write $\bW_I(t)$ for
the restriction to sites $I$ of $\bW(t)$.

\subsection{Qualitative properties of the wave process}
\label{sec-CSP}
The first two properties are intuitively obvious from the definitions;
we will not spell out details.

{\em Consistency.}
One could define the wave process on a finite site space 
$\{0,1,2,\ldots,i_0\}$.
These processes are consistent as $i_0$ varies.
In other words, taking first the wave process $\bW$ on sites
$\Zbold^+$, 
the restricted process $\bW_{[0,i_0]}$ is distributed as the wave process on sites
$\{0,1,2,\ldots,i_0\}$.

{\em Conditional independence.}
In the setting above, condition on the entire restricted process
$(\bW_{[0,i_0]}(t), 0 \leq t < \infty)$
and on the times $\sigma_j$ and right endpoints $r_j$
of waves emanating from sites in $[0,i_0]$ which have
$r_j > i_0$.
Then conditionally, the process
$\bW_{[i_0+1,\infty)}(t)$
evolves as the wave process, except that at times $\sigma_j$
the sites in $[i_0+1,r_j]$ are reset to i.i.d. Bernoulli($p$).

We next spell out a slightly more subtle
{\em conditional stationarity} property.
Consider the wave process started at time $0$ with an arbitrary
initial configuration $\bx(0)$,
and with an arbitrary initial particle distinguished.
At each subsequent time
exactly one particle is distinguished, according to the following rule.
Once distinguished,
a particle (at site $j$ say) remains distinguished until
the first time some wave emanates from some site $i \in \{j-v,j_v+1,
\ldots, j-1\}$ and at that time the particle at site $i$ becomes
the distinguished particle.

Let $q(t)$ be the position of the distinguished particle at time $t$.
\begin{Lemma}
\label{Lsut}
Fix an interval of sites $I = [i_0,i_1]$.
Suppose the initial random configuration $\bW(0)$ has a particle
at site $i_0$, the distinguished particle, and 
has
i.i.d. Bernoulli($p$) distribution
on $(i_0,i_1]$.
For each $t$, conditional on $\{q(t) = i^\prime\}$,
the distribution of $\bW_{[i^\prime +1,i_1]}(t)$ (that is, of the wave process at time $t$ on sites
$\{i^\prime + 1, i^\prime + 2,\ldots,i_1\}$)
is i.i.d. Bernoulli($p$).
\end{Lemma}
{\em Proof.}
Let $0 < U_1 < U_2 < \ldots$ be the times when the distinguished
particle changes.
Let $\bW(U_j-)$ and $\bW(U_j)$ be the configurations before and after
the change.
Inductively suppose that
conditional on $U_j$ and on $\{q(U_j) = i^\prime\}$,
the distribution of $\bW_{[i^\prime +1,i_1]}(U_j)$
is i.i.d. Bernoulli($p$).
And site $i^\prime$ is occupied at time $U_j$.
It follows from the basic stationarity property 
(of the wave process on $\Zbold^+$, translated to $[i^\prime,\infty)$)
that in the absence of waves crossing into $[i^\prime,\infty)$
from below, the distribution at any subsequent time $t$ 
of $\bW_{[i^\prime +1,i_1]}(t)$ will be
i.i.d. Bernoulli($p$).
By the conditional independence property, 
the distribution of $\bW_{[i^\prime +1,i_1]}(U_{j+1}-)$ is also
i.i.d. Bernoulli($p$),
and this remains true conditionally on $U_{j+1}$.
Given $q(U_{j+1}) = i^{\prime \prime} < i^\prime$,
the configuration 
$\bW_{[i^{\prime \prime}+1,i_1]}(U_{j+1})$
consists of the existing configuration
$\bW_{[r_j+1,i_1]}(U_{j+1}-)$
(which is
i.i.d. Bernoulli($p$)),
together with the new configuration on $[i^{\prime \prime}+1,r_j]$
created by the wave, which is also 
i.i.d. Bernoulli($p$),
thus making the whole configuration
$\bW_{[i^{\prime \prime}+1,i_1]}(U_{j+1})$
have
i.i.d. Bernoulli($p$) distribution.
This carries the induction forward, and establishes the lemma.
\qed

We next give a coupling construction.
Let $\bx^0$ and $\bx^1$ be configurations on $\Zbold^+$.
Say $\bx^1$ is an {\em extension} of $\bx^0$ if $\bx^0$ has only a finite number of occupied sites
and if the two configurations coincide on $[0,r]$, where $r= r(\bx^0)$
is the position of the rightmost particle of $\bx^0$.
In other words, $\bx^1$ consists of the particles in $\bx^0$ and
(perhaps) extra particles at arbitrary positions greater than $r$.
\begin{Lemma}
\label{Lext}
Suppose $\bx^1$ is an extension of $\bx^0$.
Then there exists a coupling $((\bW^0(t),\bW^1(t)), 0 \leq t < \infty)$
of the wave processes with initial configurations $\bx^0$ and $\bx^1$
such that at each time $t$,
the configurations
$\bW^1(t)$ and $\bW^0(t)$ coincide on sites 
$[0, \sup_{0 \leq s \leq t} r(\bW^0(s))]$.
\end{Lemma}
Call this the {\em extension coupling}.

{\em Proof.}
Suppose $\bx^1$ is an extension of $\bx^0$.
Then we can couple transitions of the joint process from these configurations
by specifying that when a wave emerges from a particle (site $i$, say)
of $\bx^1$ and
creates new particles, then if site $i$ contains a particle of $\bx^0$
we copy the wave and the positions of new particles in the other process,
and otherwise do nothing. This clearly maintains the ``extension" property.
Furthermore, if $\bx^0$ and $\bx^1$ coincide on sites $[0,r^*]$
for some $r^* \geq r(\bx^0)$ then the two coupled processes will always
agree on that interval.
\qed

Now let $(\bW^0(t), 0 \leq t < \infty)$
be the wave process on $\Zbold^+$ started with only one particle at site $0$.
Let $R_t = r(\bW^0(t))$
be the site of the rightmost particle of $\bW^0(t)$.
\begin{Lemma}
\label{Lcsg}
Suppose there exists $\lambda > 0$ such that for each $i_0$
\begin{equation}
 P(\sup_{0 \leq s \leq t} R_s \leq i_0) = O(e^{-\lambda t})
\mbox{ as } t \to \infty . \label{Rexp}
\end{equation}
Then the wave process on $\Zbold^+$ 
has spectral gap
at least $\lambda$.
\end{Lemma}
{\em Proof.}
Fix $i_0$ and write $I = [0,i_0]$.
Let $\bW^0(t)$ be as above,
let $\bW(t)$ be the stationary wave process
and 
let $\bW^*(t)$ be the wave process
started with arbitrary initial configuration $\bx^*$.
By Lemma \ref{Lext}
we can couple $\bW^0$ and $\bW$ such that
\[ P(\bW^0_I(t) \neq \bW_I(t)) \leq 
 P(\sup_{0 \leq s \leq t} R_s < i_0) . \] 
Similarly 
we can couple $\bW^0$ and $\bW^*$ such that
\[ P(\bW^0_I(t) \neq \bW^*_I(t)) \leq 
 P(\sup_{0 \leq s \leq t} R_s < i_0) . \] 
So, writing $|| \cdot||$ for variation distance,
\begin{equation}
 \sup_\bx ||P(\bW_I(t) \in \cdot|\bW_I(0) = \bx)
\ - \ \pi_I(\cdot)|| \leq 2
 P(\sup_{0 \leq s \leq t} R_s < i_0) \label{tvcon}
\end{equation}
where $\pi_I$ is 
i.i.d. Benoulli($p$) distribution on $I$.
In a finite-state reversible chain, 
the spectral gap equals the exponent in the asymptotic rate
of convergence to stationarity, and so (\ref{tvcon}) and (\ref{Rexp})
imply that
the restricted chain $\bW_{[0,i_0]}$ has spectral gap at least $\lambda$.
Since this bound is uniform in $i_0$ it extends by consistency to the wave
process on sites $\Zbold^+$
(see Appendix for more details).

\subsection{The supermartingale analysis}
\label{sec-mg}
Here is the main result of section \ref{sec-wave}.
\begin{Proposition}
\label{Pwave}
If $v = v(p) > \sfrac{10}{p} +2$
then the spectral gap of the wave process on $\Zbold^+$ is at least $3/10$,
for sufficiently small $p$.
\end{Proposition}
In outline, the idea is to apply Lemma \ref{Lcsg}.
To prove (\ref{Rexp}) for some given $\lambda$ it would suffice
to show that for some $\theta > 0$
\begin{equation}
\exp(\lambda t - \theta R_t)
\mbox{ is a supermartingale.}
\label{exactMG}
\end{equation}
For then
$E \exp(-\theta R_t) \leq e^{-\lambda t}$
and so
$P(R_t \leq x) \leq e^{\theta x} e^{- \lambda t}$.
Unfortunately (\ref{exactMG}) cannot be exactly true, because there are
``bad" configurations from which $R_t$ tends to decrease rather than increase
(for instance, if sites $R_t - v$ through $R_t - v + 100$ are occupied
but sites $R_t -v + 101$ through $R_t -1$ are unoccupied).
So we use a more elaborate argument which finesses bad configurations
by establishing the supermartingale property only at embedded random times
(Lemma \ref{Lsupmg}).
Incidently, one could modify the wave process by allowing the waves to
have {\em random} length $V$, and taking $V$ to 
have geometric distribution with sufficiently large mean
the supermartingale property (\ref{exactMG}) would be easy to verify.
However, using unbounded $V$ would make the comparison argument in section \ref{sec-comp}
become more complicated.

{\em Proof of Proposition \ref{Pwave}.}
Write $v_0 = \lceil 5/p \rceil$, so $v \geq 2v_0$.
We first consider the wave process on $\Zbold^+$, started
at time $0$ with some arbitrary finite configuration of particles.
Write $R_0$ for the position of the rightmost particle at time $0$.
Let $U$ be the first time that either a wave emerges from the
particle at $R_0$ or the particle at $R_0$ is removed by a wave emanating
from another particle.
Let $T \geq U$ be the time of the first wave whose right end is 
$ \geq R_0 + v_0$.
Note that $T \neq U$ only in the case where the particle at $R_0$
is removed by a wave whose rightmost limit is between $R_0$ and $R_0 + v_0$;
this is the case whose analysis is more difficult.
Fix some $1 > \lambda > 0$ and
define
\[ M_t = \exp(\lambda t - \sfrac{p}{2}R_t).\]
\begin{Lemma}
\label{Lsupmg}
Let $p$ be sufficiently small.
For any finite initial configuration,
\[ E M_{T \wedge (U+1)} \leq \frac{e^{\lambda -1}}{1 - \lambda}  M_0 , \quad 0<\lambda<1 . \]
\end{Lemma}
Here 
$T \wedge (U+1) := \min(T,U+1)$.
The proof has three parts.
At (\ref{100}) we show
that $T\leq U+1$ is the likely alternative.
At (\ref{TU1}) we show
that $R_T$ tends to be larger than $R_0$ on $\{T \leq U+1\}$.
At (\ref{TU2}) we show
that $R_T$ will not be much larger than $R_0$ on $\{T > U+1\}$.

{\em Proof of Lemma \ref{Lsupmg}.}
We first argue
\begin{equation}
E \left(\exp(- \sfrac{p}{2}(R_T - R_0)) | T,U \right)
\leq 2 e^{-5/2} + o(1) \mbox{ on } \{T \leq U+1\} 
\label{TU1}
\end{equation}
where $o(1)$ denotes a constant tending to $0$ as $p \downarrow 0$.
At time $T$ there is some wave 
with right endpoint $y^\prime$ and
$R_T$ is stochastically larger than
$y^\prime +1- \eta$,
where $\eta$ has geometric($p$) distribution and is independent of $T,U$.
Since $y^\prime \geq R_0 + v_0$ we see that on $\{T<\infty\}$
\begin{eqnarray*}
 E \left(\exp(- \sfrac{p}{2} (R_T - R_0)) |T,U \right)
&\leq& E \exp(- \sfrac{p}{2} (v_0 +1 - \eta))\\
&\leq& e^{-5/2} E \exp(- p \eta /2)\\
&=& e^{-5/2} (2+ o(1)). 
\end{eqnarray*}
We next shall argue that for sufficiently small $p$
\begin{equation}
P(T>U+1) \leq 1/100.
\label{100}
\end{equation}
If $R_0 \leq v_0$ then $P(T \neq U) = 0$ by the note above Lemma \ref{Lsupmg},
because the rightmost end of any wave emanating from below $R_0$ must be
$\geq v \geq R_0 + v_0$.
So assume $R_0 > v_0$.
Condition on the restriction to sites $[0,R_0-v_0]$
of the process 
$ (\bW(t), 0 \leq t < \infty)  $,
and on the times 
and rightmost ends $(t_*,x_*)$ of waves ending at positions to the
right of $R_0-v_0$.
(This conditioning is denoted by ``waves" in (\ref{wavesup}) below).
In order that $T > U+1$ there must 

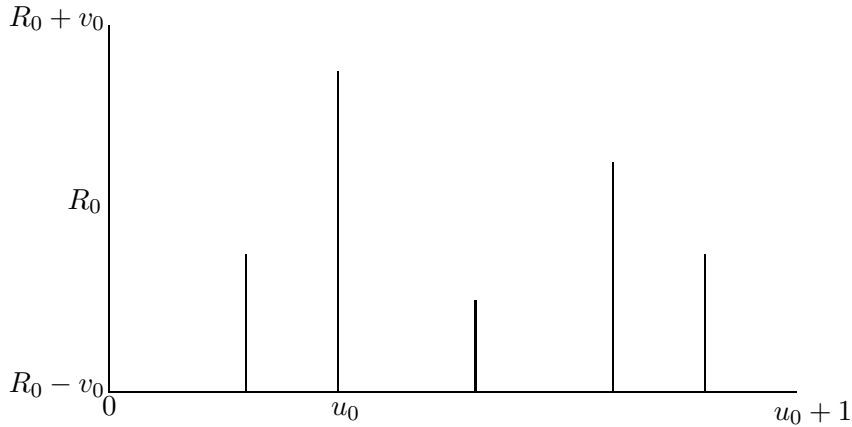
\setlength{\unitlength}{0.24in}
\begin{picture}(16,10.2)(-1,-1)
\put(0,0){\line(1,0){15}}
\put(0,0){\line(0,1){8}}
\put(-0.14,-0.5){0}
\put(-2.18,0){$R_0-v_0$}
\put(-0.9,4){$R_0$}
\put(-2.18,8){$R_0+v_0$}
\put(5,0){\line(0,1){7}}
\put(8,0){\line(0,1){2}}
\put(11,0){\line(0,1){5}}
\put(13,0){\line(0,1){3}}
\put(4.85,-0.5){$u_0$}
\put(14.5,-0.6){$u_0+1$}
\put(3,0){\line(0,1){3}}
\end{picture}

\centerline{Fig. 3. 
{\small Waves.}}

\vspace{0.12in}
\noindent
be some first wave,
at time $u_0$, with rightmost end $x_0 \in [R_0,R_0+v_0)$,
and all the other rightmost ends of waves at times before $u_0+1$ must also be
$< R_0+v_0$.  
Figure 3 illustrates such waves.

By the conditional independence property (section \ref{sec-CSP})
the restriction of the wave process to the sites $I:= (R_0-v_0,R_0+v_0]$
evolves according to the usual wave process rules, except that the
waves (on which we are conditioning) reset the sites they cover to
i.i.d. Bernoulli($p$).
The conditional probability that $T>U+1$ is the probability that no particle
in sites $I$ (neither an initial particle nor a particle created by 
any of the conditioning waves)
has a wave emanating from it before time $u_0+1$
(because any such wave would have rightmost end greater than $R_0-v_0+v \geq
R_0+v_0$).
We can upper bound this conditional probability by considering only particles
created by the wave at time $u_0$ and subsequent waves, and only sites 
$I^*:= (R_0-v_0,R_0]$.
For each site $i \in I^*$ we can decompose the time interval $[u_0,u_0+1]$
into subintervals $J$ starting at the successive times when a wave meets
site $i$.
let $\JJ$ denote the set of such intervals $J$ obtained by varying $i$,
and let $(\beta_J, J \in \JJ)$ be the i.i.d. Bernoulli($p$) random variables
indicating whether the wave created a particle at site $i$.
The number of waves created by all these particles is Poisson with conditional 
mean
$M:= \sum_J |J| \beta_J$,
where $|J|$ is the length of interval $J$.
So
\begin{equation} P(T>U+1| \mbox{ waves}) \leq E (\exp(-M) |\JJ). \label{wavesup}
\end{equation}
Now the constraints on the family $\JJ$ are
\[ |J| \leq 1 \ \forall J; \quad \sum_J |J| = v_0. \]
A routine convexity argument shows that, subject to these constraints,
$E (\exp(-M)|\JJ)$ is maximized in the case where each $|J| = 1$.  
So \[
E( \exp(-M)|\JJ) \leq E \exp(-B(v_0,p)) \]
where $B(\cdot,\cdot)$ has Binomial distribution.
As $p \to 0$ we have $v_0 p \to 5$ and so
\[ E \exp(-B(v_0,p)) \to 
\exp(- 5 (1-e^{-5})) < 1/100 \]
giving (\ref{100}).

We digress to record an elementary calculation
\begin{Lemma}
\label{Lelem}
If $Y$ is stochastically smaller than exponential(1)
then for any event $D$
\[ E \exp(Y/2)1_D \leq 2 \sqrt{P(D)} . \]
\end{Lemma}
At time $U$ there is some wave $[x,x+v]$ with $x+v \geq R_0$.
By taking time $U$ and interval $[x,x+v]$ as the initial time and
interval in Lemma \ref{Lsut}, the conclusion of that lemma easily implies
that
$R_0 - R_{U+1} + 1$ is stochastically smaller than 
geometric($p$).
This implies $p(R_0 - R_{U+1})$ is stochastically smaller than
exponential($1$).
So by applying Lemma \ref{Lelem} to $Y:= p(R_0 - R_{U+1})$
and $D:= \{T>U+1\}$ gives,
using (\ref{100}),
\begin{equation}
E \left(\exp(- \sfrac{p}{2} (R_{U+1} - R_0)) 1_{(T > U+1)} | U \right)
\leq 1/5.
\label{TU2}
\end{equation}
Now split the quantity under study in Lemma \ref{Lsupmg} over the events
$\{T \leq U+1\}$ and $\{T>U+1\}$:
\[
e^{\sfrac{p}{2} R_0} \ E M_{T \wedge (U+1)} =
E e^{\lambda T} \exp(-\sfrac{p}{2}(R_T - R_0)) 1_{(T \leq U+1)}
+ E e^{\lambda(U+1)} \exp(-\sfrac{p}{2}(R_{U+1} -R_0))
. \]
Consider the first term.
By conditioning on $T,U$ and using (\ref{TU1})
\[ \mbox{(first term)}
\leq (2e^{-5/2}+o(1)) E e^{\lambda T}1_{(T \leq U+1)}
\leq (2e^{-5/2}+o(1)) E e^{\lambda (U+1)} . \]
Similarly, by conditioning on $U$ and using (\ref{TU2}),
\[ \mbox{(second term)}
\leq \sfrac{1}{5} E e^{\lambda (U+1)} . \]
Combining these two bounds,
\[
e^{\sfrac{p}{2} R_0} \ E M_{T \wedge (U+1)}
\leq  (2e^{-5/2} + \sfrac{1}{5} + o(1)) e^\lambda 
E e^{\lambda U} . \]
The first term works out numerically to be $< e^{-1}$.
And $U$ is stochastically smaller than the exponential(1) time
at which a wave would emanate from the initial particle at $R_0$,
so $E e^{\lambda U} \leq 1/(1-\lambda)$.
This establishes Lemma \ref{Lsupmg}.
\qed

Returning to the proof of Proposition \ref{Pwave},
choose $\lambda = 3/7$ to make 
$\frac{e^{\lambda -1}}{1 - \lambda} < 1$.
Consider $(\bW^0(t), 0 \leq t < \infty)$,
the wave process on $\Zbold^+$ started with only one particle at site $0$.
Define stopping times
$0=S_0<S_1<S_2< \ldots$ by:

\begin{quote}
$S_{k+1} - S_k$ is the time $T\wedge (U+1)$ defined above
Lemma \ref{Lsupmg}, applied to the wave process 
$(\bW^0(S_k +t), \ 0 \leq t < \infty)$.
\end{quote}
So Lemma \ref{Lsupmg} implies that
$(M_{S_k}, 0 \leq k < \infty)$
is a supermartingale.
Fix $t$ and define
\[ \kappa := \min \{k: S_k \geq at\} \]
where the constant $0<a<1$ will be specified later.
Now
\begin{eqnarray*}
1 &\geq& EM_{S_\kappa} \quad \mbox{(optional sampling theorem)}\\
&\geq& EM_{S_\kappa} 1_{(R(S_\kappa) \leq x)}\\
&\geq& \exp(\lambda at - px/2) P(R(S_\kappa) \leq x) 
\end{eqnarray*}
implying
\begin{equation}
P(R_{S_\kappa} \leq x) \leq e^{px/2} e^{-\lambda at} .
\label{RSk}
\end{equation}
We next need to bound the ``overshoot"
$S_\kappa - at$.
Lemma \ref{Lovershoot} below implies
$P(S_\kappa - at \geq 1+u) \leq e^{-u}, \ 0 \leq u < \infty$
and hence
\begin{equation}
P(S_\kappa > t) \leq e \cdot e^{-(1-a)t} . \label{Skt}
\end{equation}
Finally,
\begin{eqnarray*}
P\left( \sup_{0\leq s \leq t} R_s \leq x \right) &\leq&
P(R_{S_\kappa} \leq x) + P(S_\kappa > t)\\
&\leq& 2 e^{1 \wedge px/2} e^{- \max(\lambda a,(1-a)t)} \mbox{ by }
(\ref{RSk},\ref{Skt}).
\end{eqnarray*}
Having specified $\lambda = 3/7$ we now specify $a = 7/10$
and the bound is $O(e^{-3t/10})$ as $t \to \infty$.
So Proposition \ref{Pwave} follows from Lemma \ref{Lcsg}.

\begin{Lemma}
\label{Lovershoot}
$S_\kappa - at $ is stochastically smaller than $ 1+ \eta$,
where $\eta$ denotes an exponential(1) r.v.
\end{Lemma}
{\em Proof.}
Write $t_0 = at$, so that
$\kappa = \min \{k: S_k \geq t_0\}$.
Regard time $t_0$ as the present, and condition on the past
$(\bW^0(t), 0 \leq t \leq t_0)$.
The conditioning tells us the time $S_{\kappa -1}$ and the position
of the rightmost particle $R_0 =r(S_{\kappa -1})$ at that time.
Resetting time to restart at time $S_{\kappa -1}$, we are in the setting
of figure 3, with initial configuration $\bW^0(S_{\kappa -1})$.
Now $S_\kappa - S_{\kappa -1} = T \wedge (U+1)$;
we see from the conditioning that
$T \wedge (U+1) > t_1:= t_0 - S_{\kappa -1}$
(the present time is now $t_1$)
and we are interested in the distribution of the overshoot
$\zeta:= (T \wedge (U+1)) - t_1$.
If $U \leq t_1$ then $\zeta \leq 1$ and we are done, so suppose $U > t_1$.
In that case the rightmost particle at time $t_1$ is still the particle
at site $R_0$,
and the future waiting time until $U$ is at most the exponential time $\eta$
until a wave would emanate from the particle at $R_0$.
So $\zeta \leq 1 + \eta$ and we are done.
\qed 

{\em Remark.}
Lemma \ref{Lovershoot} is slightly subtle; the fact that the conditional
distributions of $S_k - S_{k-1}$ are stochastically smaller than
$1 + \eta$ is not enough to get a bound on overshoots.

\subsection{The wave process on $\Zbold$}
\label{sec-ppzz}
Going from
the wave process on $\Zbold^+$
to the wave process on $\Zbold$
involves some easy arguments which we shall just outline.
Consider first a branching random walk ${\bf B}(t)$,
in which each site may have more than one particle,
and particles independently at rate $1$ create a wave
of offspring, one at each of the $v$ sites to the right of the parent
particle, with particles never being killed.
It is well known that (in discrete time, starting from a
single particle) the position $r({\bf B}(t))$ of the rightmost
occupied site grows asymptotically at a finite linear rate
\cite{big77}, and the same arguments give the essentially weaker 
conclusion of the next lemma.
\begin{Lemma}
\label{L24}
For the branching random walk ${\bf B}(t)$ where initially
all sites in $(-\infty,0)$ are occupied by a single particle,
\[ \lim_{L \to \infty} P(r({\bf B}(t)) > L) = 0, \quad t \mbox{ fixed}. \]
\end{Lemma}
Next, observe there is a ``basic coupling" of two versions of
the wave process, as follows.
If a site $i$ is occupied in each version, make the wave-times
from $i$ and subsequent replacements be identical; for sites
which are unmatched (occupied in one version only) let the
waves occur independently.
It is easy to check the following.
\begin{Lemma}
\label{L23}
Given two initial configurations for the wave process $\bW^{(1)}(0)$ and $\bW^{(2)}(0)$,
let ${\bf B}(0)$ be the set of unmatched sites.
Then the basic coupling 
$(\bW^1(t),\bW^2(t))$
can be constructed jointly with the branching random walk
${\bf B}(t)$ such that for each $t$ the set of unmatched sites
$(\bW^1(t),\bW^2(t))$
is a subset of the set of occupied sites in ${\bf B}(t)$.
\end{Lemma}
Now write $\bW(t)$ for the stationary wave process on sites $\Zbold^+$.
Let $\sigma_L$ be the shift map taking
$(x_i, i \geq 0)$ to $(x_{i-L}, i \geq 0)$.
\begin{Lemma}
\label{Lzz}
As $L \to \infty$ the processes $\sigma_L(\bW(t))$
converge weakly to a process $\widetilde{\bW}(t)$,
the stationary wave process on sites $\Zbold$.
\end{Lemma}
{\em Proof.}
It is enough to show that we can couple 
$\sigma_{L_1}(\bW(t))$
and
$\sigma_{L_2}(\bW(t))$
such that, for fixed $t$ and $i_0$, as
\[ L_1 \to \infty, \ L_2 \to \infty, \ L_1<L_2 \]
we have
\begin{equation}
 P( \left. \sigma_{L_1}(\bW(t)) \right|_i
= \left. \sigma_{L_2}(\bW(t)) \right|_i
\ \forall i \geq i_0) \to 1 .
\label{sss}
\end{equation}
Use the basic coupling above, where the initial configurations
coincide except on sites $[-L_2,-L_1)$, and Lemma \ref{L23}
to show that the probability in (\ref{sss}) is at least
\[ P(r({\bf B}(t)) <  i_0)  \]
where ${\bf B}(0)$ has no particles outside $[-L_2,-L_1)$.
Use Lemma \ref{L24} to show this probability $\to 1$ as $L_1 \to \infty$.
\qed

 From the ``maximal correlation" interpretation (\ref{max-cor})
of relaxation time, using functionals depending on only finitely many 
sites, one sees that for a weakly convergent sequence of interacting
particle systems, the relaxation time of the limit is at most
the limit of the relaxation times.
So the bound from Proposition \ref{Pwave}
goes through to the limit in Lemma \ref{Lzz}:
\begin{Corollary}
\label{Cz}
If $v = v(p) > \sfrac{10}{p} +2$
then the spectral gap of the wave process on $\Zbold$ is at least $3/10$,
for sufficiently small $p$.
\end{Corollary}

\section{The comparison argument}
\label{sec-comp}
\subsection{The general inequality}
Proposition \ref{Pcomp} states the general inequality we use,
in our setting of the East process $\bX$ and the wave process
$\bW$ (on state space $\Zbold$, with the same parameter $p$).

Let $Q^{\bX}(\cdot)$ be the equilibrium flow measure on 
the space $\XX$ of possible transitions 
$(\bx,\bx^\prime)$ 
of $\bX$.
That is, the first marginal of $Q^{\bX}$ is the stationary,
Bernoulli($p$), law and the conditional law is the transition
rate.
Let $Q^{\bW}(\cdot)$ be the equilibrium flow measure for $\bW$. 
For each possible transition 
$(\bw,\bw^\prime)$ of $\bW$,
define a path
$\bw = \bx(0),\bx(1),\ldots,\bx(l) = \bw^\prime$
whose steps $(\bx(i-1),\bx(i))$ are possible transitions of $\bX$,
and write $N_{(\bw,\bw^\prime)}(\cdot)$ 
for the counting measure 
on $\XX$ which counts the transitions $(\bx(i-1),\bx(i))$.
Then define
a measure on $\XX$ by:
\[ \widetilde{Q}(\cdot) := 
\int 
N_{(\bw,\bw^\prime)}(\cdot) 
\ Q^{\bW}(d \bw, d \bw^\prime) .\]
\begin{Proposition}
\label{Pcomp}
Suppose we can choose paths such that, for constants
$B, L < \infty$,\\
(i) each path length is at most $L$\\
(ii) the density
$d\widetilde{Q}/dQ^{\bX}$ is bounded a.e. by $B$.\\
Then the relaxation times of the two processes satisfy
\[ \tau(\bX) \leq BL \ \tau(\bW) . \]
\end{Proposition}
Diaconis and Saloff-Coste \cite{DSC93a} 
Theorem 2.1
prove this in the discrete time, finite state space setting,
but since the argument rests only on the general 
variational characterization (\ref{EC}) and a Cauchy-Schwarz bound,
it extends to our setting without essential alteration.
See Lemmas 1.13 -- 1.17 of Holley \cite{holley85}.

\subsection{Applying the comparison inequality}
The ``distinguished paths" required to implement the
comparison method
are readily constructed in terms of the ``minimum energy"
paths of Lemma \ref{Lm}.
At many places the bounds are crude.
\begin{Lemma}
\label{Lm2}
Let $m \geq 1$.
Take site space $[0,2^m]$ with site $0$ always occupied.
Let $\bW$ and $\bW^\prime$ be independent with 
Bernoulli($p$) distribution on sites $[1,2^m]$.
Then we can construct a path from
$\bW$ to $\bW^\prime$ of length $2 \cdot 3^m$, using only possible transitions of the
East process,
such that for each configuration $\bx$
\begin{equation}
E(\mbox{number of exits of path from $\bx$})
\leq 2 \cdot 3^m \ p^{|\bx|-m-2} . \label{Gbd}
\end{equation}
\end{Lemma}
{\em Proof.}
Write $(\hat{\bx}(u), 0 \leq u \leq 3^m)$
for the 
path from
configuration $\delta_0$
to configuration $\delta_0 + \delta_{2^m}$
given in Lemma \ref{Lm}.
Given an arbitrary configuration $\bw$ with $0$ occupied, 
for each $i \in [0,2^m]$ 
write
$s_i = \min \{u \geq 0: \hat{x}_i(u) = 1\}$ and set
\[ x_i(u) = \hat{x}_i(u) + 1_{(w_i=1)}1_{(u<s_i)} . \]
As illustrated in figure 4,
this constructs a path $(\bx(u))$ of length $3^m$ from $\bw$ to
$\delta_0 + \delta_{2^m}$;
sites initially occupied in $\bw$ remain occupied until
the $\hat{\bx}$ path first makes the site occupied, but then
behave as in the $\hat{\bx}$ path.
Joining two such paths back-to-back constructs a path
of length $2 \cdot 3^m$ between arbitrary configurations
$\bw$ and $\bw^\prime$.
It is now enough to show that, for a fixed step $u$
(w.l.o.g. $u \leq 3^m -1$),
when $\bw$ has the Bernoulli($p$) distribution,
\[ P(\bx(u) = \bx) \leq p^{|\bx|-m-2} \quad \forall \bx . \]
By Lemma \ref{Lm} we have $|\hat{\bx}(u)| \leq m+2$.
So if $|\bx| > m+2$ then there are at least 
$|\bx|-m-2$ sites which are occupied in $\bx$ but not in $\hat{\bx}(u)$;
in order for $\bx(u) = \bx$ 
it is necessary that all these sites be occupied in the initial $\bw$,
which has chance $p^{|\bx|-m-2}$.
\qed

\setlength{\unitlength}{0.12in}
\begin{picture}(-33,13)(-5,0)
\multiput(-5,9)(0,-1){10}{\circle*{0.4}}
\multiput(-4,8)(0,-1){2}{\circle*{0.4}}
\multiput(-4,2)(0,-1){2}{\circle*{0.4}}
\multiput(-3,7)(0,-1){6}{\circle*{0.4}}
\multiput(-2,9)(0,-1){6}{\circle*{0.4}}
\multiput(-1,4)(0,-1){5}{\circle*{0.4}}
\multiput(1,9)(0,-1){10}{\circle*{0.4}}
\multiput(10,9)(0,-1){10}{\circle*{0.4}}
\multiput(14,9)(0,-1){10}{\circle*{0.4}}
\multiput(15,8)(0,-1){2}{\circle*{0.4}}
\multiput(15,2)(0,-1){2}{\circle*{0.4}}
\multiput(16,9)(0,-1){8}{\circle*{0.4}}
\multiput(17,5)(0,-1){2}{\circle*{0.4}}
\multiput(18,4)(0,-1){5}{\circle*{0.4}}
\multiput(25,9)(0,-1){10}{\circle*{0.4}}
\multiput(33,9)(0,-1){10}{\circle*{0.4}}
\multiput(26,1)(0,1){2}{\circle*{0.4}}
\multiput(26,7)(0,1){2}{\circle*{0.4}}
\multiput(27,2)(0,1){6}{\circle*{0.4}}
\multiput(28,4)(0,1){2}{\circle*{0.4}}
\multiput(29,5)(0,1){5}{\circle*{0.4}}
\put(-5.23,10){0}
\put(-4.23,10){1}
\put(-3.23,10){2}
\put(-2.23,10){3}
\put(-1.23,10){4}
\put(-0.23,10){5}
\put(0.77,10){6}
\put(1.77,10){7}
\put(2.77,10){8}
\put(9.77,10){0}
\put(10.77,10){1}
\put(11.77,10){2}
\put(12.77,10){3}
\put(13.77,10){4}
\put(14.77,10){5}
\put(15.77,10){6}
\put(16.77,10){7}
\put(17.77,10){8}
\put(24.77,10){0}
\put(25.77,10){1}
\put(26.77,10){2}
\put(27.77,10){3}
\put(28.77,10){4}
\put(29.77,10){5}
\put(30.77,10){6}
\put(31.77,10){7}
\put(32.77,10){8}
\end{picture}

\vspace{0.12in}

Fig. 4. 
{\small 
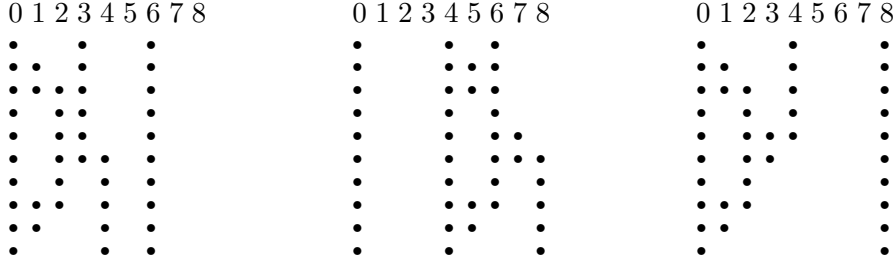
A path from $\delta_0 + \delta_3 + \delta_6$
to $\delta_0 + \delta_8$.
}

\vspace{0.1in}
\noindent
Return to the setting of Proposition \ref{Pcomp},
and suppose the wave process has waves of length $v = 2^m$.
A possible transition 
$(\bw,\bw^\prime)$ of the wave process
involves a wave from some site $i$ and only affects sites
$[i,i+2^m]$.
Use the path from $\bw$ to $\bw^\prime$ 
constructed in Lemma \ref{Lm2}.
So condition (i) in Proposition \ref{Pcomp} holds with 
\begin{equation}
L = 2 \cdot 3^m
\label{L-def}
\end{equation}
We shall argue that condition (ii) holds with
\begin{equation}
 B = 2^m \cdot 2 \cdot 3^m \ p^{-m-2} (1-p)^{-2^m} . \label{B-def}
\end{equation}
To argue this, fix some transition 
$(\bx,\bx^\prime)$ of the East process, and suppose it is site $j$
that flips in this transition.
For this to be a step along the distinguished path from
$\bw$ to $\bw^\prime$,
the wave involved in the transition $(\bw,\bw^\prime)$
must come from some site $i \in [j-2^m,j-1]$.
We get an exact expression
\begin{equation}
 \frac{d\widetilde{Q}}{dQ^{\bX}} (\bx,\bx^\prime)
= \sum_{i=j-2^m}^{j-1} 1_{(x_i = 1)}
\frac{p G(\bx_{[i,i+2^m]})}{\pi(\bx_{[i,i+2^m]}) q} 
\label{ddexact}
\end{equation}
where\\
$\bullet$ $\bx_{[i,i+2^m]}$ denotes $\bx$ restricted to sites $[i,i+2^m]$;\\
$\bullet$ $\pi(\bx_{[i,i+2^m]})$ is its Bernoulli($p$) probability;\\
$\bullet$ $q (= p \mbox{ or } 1)$ is the transition rate from $\bx$ to $\bx^\prime$;\\
$\bullet$ $G(\bx_{[i,i+2^m]})$ is the expected number of transitions
from $\bx_{[i,i+2^m]}$ to $\bx^\prime_{[i,i+2^m]}$
in the Lemma \ref{Lm2} path between two configurations which have site $i$ occupied and
which are independent Bernoulli($p$) on $[i,i+2^m]$.

Use the inequalities
\begin{eqnarray*}
q &\geq& p \\
\pi(\bx_{[i,i+2^m]}) &\geq& p^{|\bx_{[i,i+2^m]}|} (1-p)^{2^m}\\
G(\bx_{[i,i+2^m]}) &\leq& 2 \cdot 3^m \ p^{|\bx_{[i,i+2^m]}|-m-2} \mbox{ (Lemma
\ref{Lm2})}\\
1_{(x_i=1)} &\leq& 1 
\end{eqnarray*}
to show that the right side of (\ref{ddexact}) is bounded
by the quantity in (\ref{B-def}).

{\em Proof of Theorem \ref{T1}(a).}
Combining Proposition \ref{Pcomp}, with estimates (\ref{L-def},\ref{B-def}),
with Corollary \ref{Cz} shows that if $p$ is sufficiently
small and 
\[ 2^m \geq \sfrac{10}{p} +2 \]
then
\[ \tau(p) \leq \sfrac{10}{3} \cdot 2^2 \cdot 18^m 
\ p^{-m-2} (1-p)^{-2^m} . \]
Choosing the smallest $m=m(p)$ which satisfies its constraint,
\[ \tau(p) \leq \beta(1/p) \ (1/p)^m \mbox{ as } p \downarrow 0 \]
where $\beta(1/p)$ is polynomial in $1/p$.
This establishes Theorem \ref{T1}(a).

\section{Proof of the lower bound}
\subsection{Overview of argument}
\label{sec-overview}
We shall derive the lower bound by applying the variational
characterization (\ref{EC}) to a suitable test function $g$.
In many settings, some ``simple and intuitively natural" choice of test function $g$
gives a good lower bound, but here we  
proceed more indirectly by using
a function defined in terms of a different stochastic process.
See section \ref{sec-remarks} for motivation for considering
this particular process.

Fix $p$ and throughout this section set $n = \lfloor 1/p \rfloor$.
Let $S \subseteq \{1,2,\ldots, n\}$ be a non-empty set of {\em sites}.
Define the {\em coalescing random jumps} (CRJ) process with initial state $S$ as follows.
Initially there is one particle at each site in $S \cup \{0\}$.
Each particle in $S$ dies at some random time,
determined by the rule
\begin{quote}
The rate at which the particle at site $i$ dies at time $t$
equals $p^{D(t,i)}$,
where $D(t,i): = \min \{i-j: 0 \leq j < i, \mbox{ and } j \mbox{ alive at time $t$ } \} $.
\end{quote}
When particle $i$ dies, we shall say it {\em coalesces} with
the particle at the nearest lower-numbered site $j$.
Eventually all particles will coalesce with the particle at site $0$.
Let $L = L(S)$ be the random site occupied by the last-to-die particle
in $S$.  Finally, for a configuration
$\bx = (x_1,\ldots,x_n\} \in \{0,1\}^n$ define
\[ g(\bx) = P(L(S) > n/2), \mbox{ where } S = \{i: x_i = 1\}  \]
setting $g (\bx) = 0$ when $\bx$ is the zero vector.

It is very easy to estimate the required variance in the
variational characterization (\ref{EC}).
Since
\[
P(g(\bX(0)) = 0) \geq 
P(X_i(0) = 0 \ \forall 1 \leq i \leq n)
\to e^{-1} \mbox{ as } p \downarrow 0 \]
and
\[ P(g(\bX(0)) = 1) \geq
P(X_i(0) = 0, \ 1 \leq i \leq n/2; \ X_j(0) = 1,
\mbox{ for some } n/2 < j \leq n)
\] \[
\to e^{-1/2}(1-e^{-1/2}) \mbox{ as } p \downarrow 0 \]
we see
\begin{equation}
 \liminf_{p \downarrow 0} \  {\rm var\ } g(\bX(0)) > 0 . \label{var-bound}
\end{equation}
So the issue is to upper bound
${\cal E}(g,g)$.
This is done in the following three lemmas, whose proofs
are deferred.
In brief, the idea is to define a notion of ``good"
configurations, and estimate separately the contributions to
${\cal E}(g,g)$ from transitions involving good and not-good
configurations.

Write $a = \lceil 4 \log 1/p \rceil$.
Call a pair $(S,i)$ {\em admissible} if $i \in S \subseteq \{1,2,\ldots,n\}$ and 
either $i-1 \in S$ or $i = 1$.
Call an admissible pair {\em good} if there exist $k_1, k_2 \in S$ such that\\
(i) the interval $[k_1,k_2]$ contains $i$ and (if $i \neq 1$)
contains $i-1$.\\
(ii) If $k_2 < 2k_1 -a$ then $S$ does not intersect $[k_1 - b,k_1) \cup (k_2,k_2 +b]$,
where $b = k_2-k_1+a < k_1$.\\
(iii) If $k_2 \geq 2k_1 -a$ then $k_2 \leq \sfrac{n-a}{2}$
and $S$ does not intersect $(k_2,2k_2 + a]$.

Let ${\bf S}_i$ be the random subset of $\{1,2,\ldots,n\}$ containing site $i$, and (if $i\neq 1$)
containing site $i-1$, and where each other site $j$ is in ${\bf S}_i$ with 
probability $p$, independently as $j$ varies.  
So $({\bf S}_i,i)$ is admissible, by construction.
\begin{Lemma}
\label{L1}
\[ P(({\bf S}_i,i) \mbox{  is not good }) \leq \alpha(p), \]
where
$\alpha(p) > 0$ is a function satisfying
\[ \log \alpha(p) \sim - \sfrac{\log^2 (1/p)}{2 \log 2}
\mbox{ as } p \downarrow 0. \]
\end{Lemma}

Next consider two CRJ processes started from two different
initial configurations.
One can couple (i.e. define jointly) the two processes
so that, whenever a particle at site $i$ coalesces
at some time $t^\prime$
with a particle at site $j$ in one process, 
if sites $i$ and $j$ are occupied and the
intervening sites unoccupied in the other process,
then at the same time $t^\prime$ the particle at $i$ 
coalesces with the particle at $j$ in the other process.
The remaining details of the coupling are unimportant.

\begin{Lemma}
\label{L2}
Let $(S,i)$ be admissible and good, and let $|S| \geq 2$.
In the coupling
of the CRJ processes started from 
$S$ and from $S \setminus \{i\}$,
\[ P(L(S) \neq L(S \setminus \{i\})) \leq
\beta(p)   \]
where
$\beta(p) > 0$ is a function satisfying
\[ \log \beta(p) \sim - 2 \log^2 (1/p)  
\mbox{ as } p \downarrow 0. \]
\end{Lemma}

\begin{Lemma}
\label{L3}
$ {\cal E}(g,g) \leq 
\alpha(p) + \beta(p)  $.
\end{Lemma}
Lemma \ref{L3} and
the $p \downarrow 0$ asymptotics in Lemmas \ref{L1} and \ref{L2},
combined with the variational characterization (\ref{EC}) and
the variance bound (\ref{var-bound}),
establish Theorem \ref{T1}(b).

\subsection{Proofs of the lemmas}
\label{sec-pfL}
{\bf Proof of Lemma \ref{L3}.}
Consider the restriction to sites $\{1,2,\ldots,n\}$
of the East process.
Write $\pi$ for the Bernoulli($p$) stationary distribution.
The possible ``upwards" transitions are exactly the transitions
$S \setminus \{i\} \to S$ for admissible $(S,i)$.
Here we identify a configuration $\bx \in \{0,1\}^n$
with the subset $S_{\bx}:= \{i: x_i =1\}$.
The stationary flow rate for such a transition is exactly
$p \pi(S \setminus \{i\})$
if $i \geq 1$, and at most this quantity if $i=1$.
Since the contributions to ${\cal E}(g,g)$ from
upwards and downwards transitions are equal,
\[ {\cal E}(g,g) \leq 
\sum_S \sum_{i:(S,i) {\rm \ admissible}}
p \pi(S \setminus \{i\})
\ (g(S) - g(S \setminus \{i\}))^2 . \]
Recall the definition of ${\bf S}_i$ in Lemma \ref{L1}.
Since $P({\bf S}_i = S) = \pi(S)/p^2$ for admissible $(S,i)$,
and $\pi(S) = \pi(S \setminus \{i\}) p/(1-p)$,
we find
\[ {\cal E}(g,g) \leq 
p^2(1-p)
\sum_S \sum_{i:(S,i) {\rm \ admissible}}
P({\bf S}_i = S)
\ (g(S) - g(S \setminus \{i\}))^2 . \]
Bounding $\sum_i ( \cdot)$ by $n \max_i (\cdot)$
and observing $np^2(1-p) < 1$,
\[ {\cal E}(g,g) \leq 
\max_i \ E \left[ ( g({\bf S}_i) - g({\bf S}_i \setminus \{i\}) )^2
\ 1_{(({\bf S}_i,i) {\rm \ admissible})} \right] . \]
Recall $0 \leq g \leq 1$.
By Lemma \ref{L1} the contribution to the expectation
from the event where $({\bf S}_i,i)$ is {\em not} good is at most $\alpha(p)$.
By Lemma \ref{L2} the contribution to the expectation
from the event where $({\bf S}_i,i)$ {\em is} good is at most $\beta(p)$.
So
\[ {\cal E}(g,g) \leq 
\alpha(p) + \beta(p)  \]
as required.
Note that the requirement $|S| \geq 2$ in Lemma \ref{L2}
eliminates only the case $(S,i) = (\{1\},1)$,
which makes zero contribution to ${\cal E}(g,g)$.

\vspace{0.15in}

{\bf Proof of Lemma \ref{L2}.}
Let $k_1,k_2$ be as in the definition of {\em good}.
Consider first case (ii), where $k_2 < 2k_1 - a$.
Define $t_0$ by
\[ t_0p^b = p^{a/2} ; \quad b := k_2-k_1+a . \]
Consider the CRJ process started from $S$.
The distance from $k_1$ to the nearest particle to the left
(i.e. at a lower-numbered site)
is at least $b$, so
\[ P(\mbox{ particle at $k_1$ dies before time $t_0$}) 
\leq t_0p^b = p^{a/2} . \]
Similarly, the distance from $k_2$ to the nearest particle to
the right is at least $b$, so
\[ P(\mbox{ some particle starting to the right coalesces with any } \]
\[ \mbox{particle starting in $[k_1,k_2]$ before time $t_0$}) 
\quad \leq t_0p^b = p^{a/2} . \]
Now assume that the particle at site $k_1$ has not died before time $t_0$.
Then the chance that the particle initially at site $k_2$ has {\em not}
coalesced with the particle at site $k_1$ by time $t_0$ is at most
\begin{equation}
 (k_2-k_1) \exp\left(-t_0p^{k_2-k_1}/(k_2-k_1)\right) . \label{gap}
\end{equation}
To argue (\ref{gap}), divide the time interval $[0,t_0]$
into $k_2-k_1$ subintervals of length $t_0/(k_2-k_1)$.
In order for the event in question to occur, in one of these
subintervals the particle initially at site $k_2$ must not
coalesce with any other particle.
But the coalescence rate is at least $p^{k_2-k_1}$, so the
chance of non-coalescence over a subinterval is at most
$\exp(-p^{k_2-k_1} t_0/(k_2-k_1))$.

Since $k_2 - k_1 < n \leq 1/p$,
the quantity in (\ref{gap}) is at most
$p^{-1} \exp(-p^{1-a})$.
Combining these estimates, the chance that the event
\begin{quote}
(*) at time $t_0$ all the particles initially in $[k_1,k_2]$,
and no other particles, have coalesced into a particle 
currently at site $k_1$
\end{quote}
{\em fails}
is at most
$\beta(p) /2$, where
$\beta(p):= 4p^{a/2} + 
2p^{-1} \exp(-p^{1-a})$.
Now the same argument gives the same bound for the 
CRJ process started from $S \setminus \{i\}$,
replacing $k_2$ by $i-1$ in the case $k_2=i$.
It follows that, in the coupling, 
outside an event of probability 
$\beta(p)$
the two processes are equal at time $t_0$,
implying that $L(S) = L(S \setminus \{i\})$.
This establishes Lemma \ref{L2} in the case
(ii).

Case (iii) is similar.
Take $t_0 = p^{k_2+a/2}$ and consider the event
\begin{quote}
(**) at time $t_0$ all the particles initially in $[1,k_2]$,
and no other particles, have coalesced into a particle 
currently at site $0$.
\end{quote}
Arguing as above, the chance that event (**) fails 
can be bounded by
$p^{a/2} +
p^{-1} \exp(-p^{1-a})$.
Here the first term bounds the chance that some particle initially
to the right of $k_2$
(and hence, by (iii), to the right of $2k_2+a$)
coalesces with the particle initially at $k_2$;
the second term bounds the chance that the particle initially at
$k_2$ has not coalesced with the particle at site $0$.
The remainder of the argument follows the previous case.

\vspace{0.15in}

{\bf Proof of Lemma \ref{L1}.}
We start by examining some deterministic properties of an 
admissible pair $(S,i)$ which is {\em not} good.
We do the case $i>1$; the case $i=1$ is similar.
Set $[k_1(0),k_2(0)] = [i-1,i]$ and inductively for $m = 0,1,2,3,\ldots$
specify

(a)  if $k_2(m)<2k_1(m)-a$ then
\[ [k_1(m+1),k_2(m+1)] = [j,k_2(m)] \mbox{ or } [k_1(m),j] \]
where $j$ is the site in $S \setminus [k_1(m),k_2(m)]$
closest to the interval $[k_1(m),k_2(m)]$,
breaking ties arbitrarily.

(b) If $2k_1(m)-a \leq k_2(m) \leq \sfrac{n-a}{2}$ then
\[ [k_1(m+1),k_2(m+1)] = [k_1(m),j] \]
where $j$ is the site in $S \cap (k_2(m),n]$
closest to $k_2(m)$.

The fact that $(S,i)$ is not good implies that no pair
$\{k_1(m),k_2(m)\}$ can satisfy conditions (ii,iii),
and this implies that the inductive construction makes sense
(i.e. the required $j$'s exist at each step)
and that $(k_1(m),k_2(m))$ is well-defined for all
$0 \leq m \leq m_1$, where
\[ m_1:= \min \{m: k_2(m) \geq 2k_1(m)-a
\mbox{ and } k_2(m) > \sfrac{n-a}{2} \} . \]
Further,
for $0 \leq m < m_1$:

(c) if $k_2(m)<2k_1(m)-a$ then
$k_2(m+1)-k_1(m+1) \leq a + 2(k_2(m)-k_1(m))$

(d) if $k_2(m) \geq 2k_1(m)-a$ then
$k_2(m+1) \leq 2k_2(m) + a$.

\begin{Lemma}
\label{L4}
Define
\begin{eqnarray*}
b(m) &=& k_2(m)-k_1(m) \mbox{ if } k_2(m)<2k_1(m) -a\\
&=& k_2(m) \mbox{ if } k_2(m) \geq 2k_1(m) -a .
\end{eqnarray*}
Then
\[a + b(m) \leq 2^{m+1}(a+1), \ m \leq m_1. \]
Moreover
\[ m_1 \geq m_0:= \max \{m:\ 2^m(a+1) \leq \sfrac{n-a}{2} \} . \]
\end{Lemma}
{\em Proof.}
By (c,d), the inequality
\[ b(m+1) \leq a + 2b(m) \]
holds for all $0 \leq m < m_1$ {\em except} perhaps 
for the first value of $m$ ($m^*$, say) that
$k_2(m^*+1) \geq 2k_1(m^*+1) -a$.
For this particular value
\begin{eqnarray*}
k_2(m^*+1) &=&k_1(m^*+1) + (k_2(m^*+1)-k_1(m^*+1))\\
&\leq&k_1(m^*+1) + a + 2b(m^*)\\
&\leq& \sfrac{k_2(m^*+1)+a}{2} + a + 2b(m^*)
\end{eqnarray*}
and rearranging gives
\[ b(m^*+1) = k_2(m^*+1) \leq 
a + 2(a+2b(m^*)) . \]
In other words, if $\bar{b}(m)$ solves the recursion
\[ \bar{b}(m+1) = a + 2\bar{b}(m); \quad \bar{b}(0) = 1 \]
then inductively
\begin{eqnarray*}
b(m) &\leq& \bar{b}(m), \ m \leq m^*\\
b(m^*+1)&\leq& \bar{b}(m^*+2) \\
b(m) &\leq& \bar{b}(m+1), \ m^* < m \leq m_1 .
\end{eqnarray*}
But explicitly
\[ \bar{b}(m) = 2^m + a(2^m-1) \] 
and so
\[ a + b(m) \leq a + \bar{b}(m+1) = 2^{m+1}(1+a) \]
establishing the first inequality in the lemma.
For the second inequality, from the definition of $m_1$ we have
$b(m_1) > \sfrac{n-a}{2}$ and hence
$2^{m_1+1}(1+a) > \sfrac{n-a}{2}$.
So $m_1 \geq m_0$ by definition of $m_0$.
\qed

Returning to the proof of Lemma \ref{L1},
consider the random set ${\bf S}_i$.
In order that $({\bf S}_i,i)$ be not good it is necessary that
the random $[k_1(m),k_2(m)], \ 0 \leq m \leq m_0$
constructed by (a,b) are well-defined and satisfy (c,d).
The conditional probability ($\alpha_m$, say)
this holds for $k_1(m+1),k_2(m+1)$ is at most

(case (a,c)): $2(a+k_2(m)-k_1(m))p = 2(a+b(m))p$

(case (b,d)): $(a+k_2(m))p = (a+b(m))p$.

\noindent
By Lemma \ref{L4} $a+b(m) \leq 2^{m+1}(1+a)$
and so
$\alpha_m \leq 2^{m+2}(1+a)p$.
So the unconditional probability that $({\bf S}_i,i)$ is not good
is at most
\begin{eqnarray*}
\alpha(p) &:=& \prod_{m=0}^{m_0-1} 2^{m+2}(a+1)p\\
&\leq&  2^{(m_0+1)(m_0+2)/2} ((a+1)p)^{m_0} .
\end{eqnarray*}
 From the definition of $m_0$ we have
\[ 2^{m_0+1} \leq \sfrac{n}{a+1} \leq \sfrac{1}{p(a+1)} \]
leading to
\[ \alpha(p) \leq  \ \left( p(a+1)\right)^{m_0/2 -1} . \]
Then as $p \downarrow 0$
\[ \log \alpha(p) \sim \sfrac{1}{2} m_0 \log p
\sim \sfrac{1}{2} \log_2 (1/p) \times \log p \]
establishing Lemma \ref{L1}. 

\subsection{Remarks on the proof of the lower bound}
\label{sec-remarks}
The CRJ process is designed as a caricature of the East process
started from all $1$'s, where the occupied sites of the CRJ
process at time $t$ mimic the sites in the East process which have
been occupied {\em throughout} the time interval $[0,t]$.
The exact details of the CRJ process seem irrelevant for our
argument.

There is a shorter argument which leads to a cruder lower bound.
Take $p = 2^{-m}$ and sites $\{0,1,2,\ldots,2^m\}$.
Apply the variational characterization (\ref{EC})
to $g = 1_A$ where $A$ is the set of configurations reachable
from the basic (only site $0$ occupied) configuration by
paths using no more than $m-2$ extra particles.
Using straightforward calculations, and combinatorial lemmas
analogous to Lemma \ref{Lm}, one can prove
that for each $a>0$
\[ \tau(p) \geq  (\sfrac{1}{p})^{a \log \log (1/p)} 
\mbox{ for all sufficiently small } p . \]
Heuristics along these lines were given in \cite{mauch}.

\vspace{0.16in}

{\em Acknowledgements.}
This work arose from discussions with Hans C. Andersen of his
work \cite{PYA00}.
Authors' research supported by N.S.F. Grants
MCS 99-70901 and 95-04379.

\newpage
\appendix
\section{Appendix: some technical background}
The East process takes values in the space
${\cal X} = \{0,1\}^\Zbold$.
Give this space the usual product topology and $\sigma$-algebra,
and let $\pi$ be product Bernoulli($p$) measure.
Write $L^2$ for the space of $\pi$-square-integrable functions
$f: {\cal X} \to R$
and write ${\cal C}$ for the space of continuous functions
$f: {\cal X} \to R$.
Finally let
\[
 {\cal D}_1 = \{f \in {\cal C}: \sum_k \sup_{\bx} 
|f(\bx^k) - f(\bx)| < \infty \}   
\]
where $\bx^k$ is the configuration obtained from $\bx$ by
flipping the $k$'th coordinate.
The space ${\cal D}_1$ includes the space ${\cal D}$ of functions
depending on only finitely many coordinates.
It will serve as a core (see e.g. Liggett \cite{lig85} Chapter 1.3) for the generator of the East process.

For $k \in \Zbold$ and $\bx \in {\cal X}$ define 
a measure on $\{0,1\}$ by
\[ c_k(\bx,y) = (1-p)1_{(y=1,x_{k-1}=1)} + p1_{(y=0,x_{k-1}=1)} . \]
Define a linear operator $\Omega$ on ${\cal D}$ by
\[
\Omega f(\bx) = \sum_{k = - \infty}^\infty 
\sum_{y=0}^1 1_{(x_k=y)}
 \left(f(\bx^{k}) - f(\bx) \right) \ c_k(\bx, y)
 . \]
Theorem 3.9 of \cite{lig85} shows
(the conditions (3.3) and (3.8) being easy to check)
that the closure $\bar{\Omega}$ of $\Omega$ is a 
generator of a Markov
semigroup $S_t$ on ${\cal C}$.
This extends (\cite{lig85}) Proposition 4.1) to a generator and semigroup on
$L^2$ with ${\cal D}_1$ as core; denote these also by
$\bar{\Omega}$ and $S_t$.
The semigroup specifies an ${\cal X}$-valued stochastic process $\bX_t$
such that
\[ E(f(\bX_t)|\bX_0 = \bx) = S_tf(\bx); \quad
t \in [0,\infty), \ f \in L^2 . \]
This is a precise construction of the East process.
The operator $\bar{\Omega}$ is a self adjoint unbounded operator on $L^2$
with ${\cal D}_1$ as a core (\cite{lig85} Chapter 4.4).
The spectral theorem implies
\[  \bar{\Omega} = \int_0^\infty \lambda G(d\lambda) \]
where $G(\lambda)$ is a resolution of the identity; that is,
a family of projections on $L^2$ satisfying
\begin{eqnarray*}
 G(\lambda_1)G(\lambda_2) &=& G(\lambda_1 \wedge \lambda_2)\\
\lim_{\lambda \to \infty} G(\lambda)f &=& f\\
\lim_{\lambda \to 0} G(\lambda)f &=& G(0)f .
\end{eqnarray*}
Writing $\sigma(\bar{\Omega})$ for the support of $G(\cdot)$,
the spectral gap is defined as 
\[{\rm gap}(\bar{\Omega}) :=\min \{\lambda > 0: \ \lambda \in \sigma(\bar{\Omega}) \} . \]

When we have weak convergence of reversible processes
(generators $\bar{\Omega}_n$ and $\bar{\Omega}$ say),
we would like to conclude that 
\begin{equation}
 {\rm gap}(\bar{\Omega}) \geq \limsup_n {\rm gap}(\bar{\Omega}_n) . \label{RS}
\end{equation}
Reed and Simon \cite{Reed-Simon-1}
show that (\ref{RS}) holds provided the generators 
$\bar{\Omega_n}, \bar{\Omega}$ have ${\cal D}$ as a common core and
provided 
\[ || \bar{\Omega}_nf - \bar{\Omega}f||_2 \to 0 \ \forall f \in {\cal D}_1 . \]
This can readily be checked for the East process on $\{0,1,2,\ldots,n\}$ and the
East process on $\{0,1,2,\ldots\}$.  
In the paper we applied (\ref{RS}) to the wave process
(at the end of sections \ref{sec-CSP} and \ref{sec-ppzz}),
but this is just a similar argument.

\end{document}